\documentclass[11pt]{amsart}

\usepackage{amsfonts,epsfig}
\usepackage{latexsym}
\usepackage{amssymb}
\usepackage{amsmath}
\usepackage{amsthm}
\usepackage{graphics}
\usepackage[all]{xy}
\usepackage[T2A]{fontenc}
\usepackage{multirow}
\usepackage{array, booktabs, ctable}
\usepackage{enumerate}
\usepackage{booktabs}
\usepackage{longtable}
\usepackage{multirow}
\usepackage[pagebackref]{hyperref}
\usepackage{xcolor}

\newtheorem{theorem}{Theorem}

\newtheorem{proposition}[theorem]{Proposition}

\newtheorem{remark}[theorem]{Remark}

\newcommand{\cC}{{\mathcal{C}}}

\newcommand{\Q}{{\mathbb{Q}}}
\newcommand{\Qbar}{\overline{\mathbb Q}}
\newcommand{\Z}{{\mathbb{Z}}}

\newcommand{\fcm}{{\mathfrak{cm}}}

\DeclareMathOperator{\cm}{CM} 
\DeclareMathOperator{\tors}{tors}

\title[Torsion of rational elliptic curves with CM over quadratic fields]{Explicit characterization of the torsion growth of\\ rational elliptic curves with complex multiplication\\ over quadratic fields}

\author{Enrique Gonz\'alez--Jim\'enez}
\address{Departamento de Matem{\'a}ticas, Universidad Aut{\'o}noma de Madrid, Madrid, Spain}
\email{enrique.gonzalez.jimenez@uam.es}
\urladdr{http://matematicas.uam.es/~enrique.gonzalez.jimenez}

\thanks{The author was partially  supported by the grant PGC2018-095392-B-I00.}

\begin{document}
%\date{\today}

\maketitle

\begin{abstract}
In a series of papers we classify the possible torsion structures of rational elliptic curves base-extended to number fields of a fixed degree. In this paper we turn our attention to the question of how the torsion of an elliptic curve with complex multiplication defined over the rationals grows over quadratic fields. We go further and we give an explicit characterization of the quadratic fields where the torsion grows in terms of some invariants attached to the curve.
\end{abstract}
%\tableofcontents

\section{Introduction}
For over a century mathematicians have been enamored with the study of elliptic curves.  Of particular interest has been  characterizing the possible torsion structures of elliptic curves defined over a number field of fixed degree. The set of possible groups (up to isomorphism) is denoted by $\Phi(d)$ and much progress in understanding this set has been made in the last few decades. Thanks to Merel \cite{Merel}, it is known that $\Phi(d)$ is finite.  Beyond this, there are only two cases published in the literature. The case when $d=1$ was proven by Mazur \cite{Mazur1978} and the case when $d=2$ by Kamienny \cite{Kamienny92} and Kenku and Momose \cite{KenkuMomose88}. {Recently, Derickx, Etropolski, van Hoeij, Morrow, and Zureick-Brown have released an article \cite{dehmz} with a complete description of $\Phi(3)$. As of now, the case when $d > 3$ remains open.}

The purpose of this paper is somewhat different: we are interested in studying how the torsion grows when the field of definition is enlarged. That is, we will restrict to the case when the field of definition of the elliptic curve is actually $\Q$, but considered over a larger number field. In this context, the first problem is to characterize the set $\Phi_\Q(d)$ of possible groups (up to isomorphism) that can appear as the torsion subgroup of an elliptic curve defined over $\Q$ base extended to a field of degree $d$.  The set $\Phi_\Q(d)$ has been completely classified for $d=2,3,4,5,7$ and for any positive integer $d$ whose prime divisors are greater than $7$ (cf. \cite{Najman16,Chou16, GJN16, GJ17}). The case $d=6$ has been studied in \cite{HGJ18,G21}. 

Another problem that we are interested in is understanding the behavior of a particular torsion group $G \in \Phi(1)$ when we enlarge the base field $\Q$.
That is, if $E/\Q$ is an elliptic curve such that $E(\Q)_{\tors}\simeq G$, what groups (up to isomorphism) are of the form $E(K)_{\tors}$ as $K$ runs over all number fields of degree $d$? 
Let  $\Phi_\Q(d,G)$ denote the subset of $\Phi_\Q(d)$ such that $E$ runs through all elliptic curves over $\Q$ such that $E(\Q)_{\tors}\simeq G$. The set $\Phi_\Q(d,G)$ has been determined for $d=2,3,4,5,7$ and for any positive integer $d$ whose prime divisors are greater than $7$ (cf. \cite{GJT14,GJNT16,GJLR18,GJN16,GJ17}). The case $d=6$ has been studied in \cite{HGJ18}.

In fact, we can refine the previous question even further. We start by noting that if $E$ is an elliptic curve defined over $\Q$ and $K$ a number field such that the torsion of $E$  grows from $\Q$ to $K$, then of course the torsion of $E$ also grows from $\Q$ to any extension of $K$. With this in mind, we say that the torsion growth over $K$ is primitive if $E(K')_{\tors}\subsetneq E(K)_{\tors}$ for any subfield $K'\subsetneq K$. We denote by $\mathcal{H}_{\Q}(d,E)$ the list formed by $E(\Q)_{\tors}$ together with the groups $H$ such that there exists a number field $K$ (up to isomorphism) of degree dividing $d$ such that $E$ has primitive torsion growth over $K$ and $E(K)_{\tors}\simeq H$. We point out here that we are allowing two (or more) of the torsion subgroups $H$ to be isomorphic if the corresponding number fields are not isomorphic. We call $\mathcal{H}_{\Q}(d,E)$ the set of torsion configurations (of degree $d$) of the elliptic curve $E/\Q$. We let $\mathcal{H}_{\Q}(d)$ denote the set of $\mathcal{H}_{\Q}(d,E)$ as $E$ runs over all elliptic curves defined over $\Q$ such that $\mathcal{H}_{\Q}(d,E)\ne \{E(\Q)_{\tors}\}$. Finally, for any $G\in \Phi(1)$ we define $\mathcal{H}_{\Q}(d,G)$ as the set of lists $\mathcal{H}_{\Q}(d,E)$ where $E$ runs over all the elliptic curves defined over $\Q$ such that $E(\Q)_{\tors}\simeq G$ and $\mathcal{H}_{\Q}(d,E)\ne \{G\}$. Denote by $h_\Q(d)$ the maximum of the cardinality of the sets $S$ when $S\in \mathcal{H}_{\Q}(d)$, in other words, $h_\Q(d)$ gives the maximum number of field extensions of degrees dividing $d$ where there is primitive torsion growth. The sets $\mathcal{H}_{\Q}(d,G)$ and $\mathcal{H}_{\Q}(d)$ and the integer $h_\Q(d)$, have been determined for $d=2,3,5,7$ and for any positive integer $d$ whose prime divisors are greater than $7$ (cf. \cite{GJT15,GJNT16,GJ17,GJN16}). The cases $d=4$ and $d=6$ have been studied in \cite{GJLR18} and \cite{HGJ18} respectively.

Finally, once we have a complete classification of the sets  $\Phi_\Q(d)$, $\Phi_\Q(d,G)$ and $\mathcal{H}_{\Q}(d,G)$ for a fixed $d$, there is the much harder problem (cf. \cite[Problem 2]{GJT14}): 

\

{\bf Problem 1:} {\em  Let $E$ be an elliptic curve defined over $\Q$ and $d$ a positive integer. Is there a precise (and easy) description of which are the possible number fields $K$ of degree dividing $d$  where the torsion growth over $K$ is primitive, ideally in terms of some invariant(s) of the curve?}

\

In  \cite{GJN19}  the authors present a fast algorithm that takes as input an elliptic curve defined over $\Q$ and an integer $d$ and returns all the number fields of degree dividing $d$ where there is primitive torsion growth. But this algorithm does not provide a solution to our problem since it does not compute the number fields in terms of the invariants of the elliptic curve.

Ideally we would like an answer to this problem in complete generality. As a first step towards that goal we describe completely the torsion growth of elliptic curves with complex multiplication (CM) over $\Q$ base-extended to a quadratic field. The case of base extending to cubic number fields is solved in \cite{GJCM3}. In an ongoing project \cite{GJCMn}, we will solve the problem for number fields of low degree ($d\le 23$).

We define  $\Phi^{\cm}(d)$, $\Phi^{\cm}_\Q(d)$, $\Phi^{\cm}_\Q(d,G)$, $\mathcal{H}^{\cm}_{\Q}(d,G)$, $\mathcal{H}^{\cm}_{\Q}(d)$, $h_\Q^{\cm}(d)$,  to be the sets and constants defined analogously to the ones above but restricted to elliptic curves with complex multiplication.

The set $\Phi^{\cm}(1)$ was determined by Olson \cite{Olson74}:
$$
\Phi^{\cm}(1)=\left\{ \cC_1\,,\,  \cC_2\,,\,  \cC_3\,,\,  \cC_4\,,\,  \cC_6\,,\, \cC_2\times\cC_2\right\}.
$$
{To the best of the author's knowledge\footnote{M\"uller, Str\"oher, and Zimmer in \cite{MSZ89}; and Fung, M\"uller, Peth\H{o}, Str\"oher, Weis, Williams, and Zimmer in \cite{FSWZ90,PWZ97} determine all torsion subgroups of elliptic curves with algebraic integer $j$-invariant over quadratic and cubic fields respectively. Note that elliptic curves with CM form a subclass of elliptic curves with integral $j$-invariant. But they do not identify the CM case within this larger classification problem.}, the first classification of the quadratic and cubic case was was done by Clark in \cite[Theorem 4]{Clark2004}. Although it appears for the first time in print in \cite{Clark2014}, where Clark, Corn, Rice, and Stankewicz computed the sets $\Phi^{\cm}(d)$, for $2\le d\le 13$.} In particular,
$$
\Phi^{\cm}(2) = \Phi^{\cm}(1)\cup \left\{ \cC_7\,,\,\cC_{10}\,,\,\cC_2\times\cC_4\,,\,\cC_2\times\cC_6\,,\,\cC_3\times\cC_3\,\right\}.
$$
{Moreover, Bourdon, Clark and  Stankewicz \cite{BCS17} determine $\Phi^{\cm}(p)$ for any prime $p$, and Bourdon and Pollack \cite{BP17} generalize to $\Phi^{\cm}(d)$ for all odd $d$, showing the answer explicitly for all odd $d < 100$.}

\

{ The main results of this paper are the following:
}

\begin{theorem}\label{teo1} $\Phi_{\Q}^{\cm}(2)=\Phi^{\cm}\left(2\right)\setminus\left\{ \cC_7\,,\,\cC_{10}  \right\}$.
\end{theorem}

\begin{theorem}\label{teo2}
Let $G \in \Phi^{\cm}(1)$. The sets $\Phi^{\cm}_\Q \left(2,G \right)$ and $\mathcal{H}^{\cm}_{\Q}(2,G)$ appear in Table \ref{Phi2}.
\begin{table}[h]
 \caption{$\Phi^{\cm}_\Q \left(2,G \right)$ and $\mathcal{H}^{\cm}_{\Q}(2,G)$ for $G \in \Phi^{\cm}(1)$}\label{Phi2}
\begin{tabular}{|c|c|l|}
\hline
$G$ & $ \Phi^{\cm}_\Q \left(2,G \right)\setminus\{G\}$ & $\mathcal{H}^{\cm}_{\Q}(2,G)$\\
\hline\hline
\multirow{2}{*}{$\cC_1$}  & \multirow{2}{*}{$\left\{\, {\cC_3}\, \right\}$} & {$\cC_3$}     \\
\cline{3-3}
& & {$\cC_3\,,\cC_3$}   \\
\hline
\multirow{5}{*}{$\cC_2$}  & \multirow{5}{*}{$\left\{\, {\cC_4}\,,\, \cC_6\,,\, \cC_{2}\times\cC_{2}\,,\, \cC_{2}\times\cC_{6}\, \right\}$} & {$\cC_2\times\cC_2$}     \\
\cline{3-3}
& & {$\cC_2\times \cC_6$}    \\
\cline{3-3}
& &  {$\cC_2\times\cC_2\,,\cC_6$}    \\
\cline{3-3}
& &  {$\cC_2\times\cC_6\,,\cC_6$}      \\
\cline{3-3}
& &  {$\cC_2\times\cC_2\,,\cC_4\,,\cC_4$} \\
\hline
 $\cC_3$   & $\left\{\,  \cC_{3}\times\cC_{3}\,\right\}$  & {$\cC_3\times\cC_{3}$}  \\
\hline
$\cC_4$ & $\left\{\,\cC_2\times{\cC_{4}}\, \right\}$ & {$\cC_2\times \cC_4$}  \\
\hline
$\cC_6$ & $\left\{\, \cC_{2}\times\cC_{6} \, \right\}$ & {$\cC_2\times \cC_6$}    \\
\hline
\multirow{2}{*}{$\cC_2\times \cC_2$}  & \multirow{2}{*}{$\left\{ \,\cC_2 \times \cC_{4}\,\right\}$}& {$\cC_2\times\cC_4$}   \\
\cline{3-3}
& & {$\cC_2\times\cC_4\,,\cC_2\times\cC_4$}   \\
\hline
\end{tabular}
\end{table}
In particular, $h^{\cm}_\Q(2)=3$.
\end{theorem}

Finally, we give an affirmative answer to Problem 1 for the case of elliptic curves defined over $\Q$ with CM base changed to quadratic fields in terms of what we define as the CM-invariants of the curve (see \S \ref{sec:CM} for the definition).

\begin{theorem}\label{teo3}
Let $E$ be an elliptic curve defined over $\Q$ with CM. Let $k$ and $\mathfrak{cm}$ its corresponding CM-invariants. Table \ref{d2} gives an explicit description of torsion growth over quadratic fields of $E$ depending only on the integers $k$ and $\mathfrak{cm}$. The 4$^\text{th}$ column gives a list of the form $H_1,\dots, H_n$ and the 5$^\text{th}$ column gives $\sqrt{d_1},\dots, \sqrt{d_n}$ such that $E(\Q(\sqrt{d_i}))_{\tors}\simeq H_i$, for $i=1,\dots,n$.
\end{theorem}

\begin{table}[ht!]
 \caption{Explicit description of torsion growth over quadratic fields of elliptic curves defined over $\Q$ with CM}\label{d2}
\begin{tabular}{|c|c|c|l|c|}
\hline
$\mathfrak{cm}$ & $k$ such that $E=E^k_{\mathfrak{cm}}$ & $  E(\Q)_{\tors}$ & $\mathcal{H}_\Q(2,E)\setminus\{E(\Q)_{\tors}\}$ &  quadratics\\
\hline\hline
 \multirow{7}{*}{$3$} & $1$ & $\cC_6$ & $\cC_2\times \cC_6$  & $\sqrt{-3}$ \\
\cline{2-5}
& $16,-432$&   \multirow{2}{*}{$\cC_3$}  & {$\cC_3\times\cC_3$}  &  $\sqrt{-3}$ \\
\cline{2-2}\cline{4-5}
 & $r^2 \,\, (r\ne\pm 1,\pm 4)$ &   & $-$ & $-$  \\
\cline{2-5}
 &$-27$  &  \multirow{2}{*}{$\cC_2$}  & {$\cC_2\times\cC_6$}  &  $\sqrt{-3}$ \\
\cline{2-2}\cline{4-5}
 &{$r^3\,\,(r\ne 1,-3)$}  &  & {$\cC_2\times\cC_2\,,\cC_6$}  & $\sqrt{-3},\sqrt{r}$  \\
\cline{2-5}
 & $2r^3\,\,(r\ne -6,2)$&   \multirow{2}{*}{$\cC_1$}  & {$\cC_3,\,\cC_3$}  &  $\sqrt{2r},\sqrt{-6r}$ \\
\cline{2-2}\cline{4-5}
 & $\ne r^2,r^3, 2r^3$ & & $\cC_3$ & $\sqrt{k}$\\
\hline
  \multirow{3}{*}{$12$} & $1$ & $\cC_6$ & $\cC_2\times\cC_6$  & $\sqrt{3}$ \\
\cline{2-5}
 & $3$ & \multirow{2}{*}{$\cC_2$} & $\cC_2\times\cC_6$  & $\sqrt{3}$\\
\cline{2-2}\cline{4-5}
  & $\ne 1,3$ & & $\cC_2\times\cC_2,\,\cC_6$  & $\sqrt{3},\sqrt{k}$\\
\hline
 \multirow{2}{*}{$27$} & $1$ & $\cC_3$ & $-$ & $-$\\
\cline{2-5}
  &$\ne 1$  & $\cC_1$ & $\cC_3$  & $\sqrt{k}$\\
\cline{4-5}
\hline
\multirow{6}{*}{$4$} & $-1$ & \multirow{3}{*}{$\cC_2\times \cC_2$} & $\cC_2\times\cC_4,\,\cC_2\times\cC_4$  & $\sqrt{-1},\sqrt{2}$\\
\cline{2-2}\cline{4-5}
& $-4$ & & $\cC_2\times\cC_4$  & $\sqrt{2}$\\
\cline{2-2}\cline{4-5}
 & $-r^2 \,\, (r\ne\pm 1,\pm 2)$ &   & $-$ & $-$  \\
\cline{2-5}
 & $4$ & {$\cC_4$} & $\cC_2\times\cC_4$  & $\sqrt{-1}$\\
\cline{2-5}
   & $r^2\,\, (r\ne \pm 2)$  &  \multirow{2}{*}{$\cC_2$}  & $\cC_2\times\cC_2\,,\cC_4\,,\cC_4$  & $\sqrt{-1},\sqrt{2r},\sqrt{-2r}$ \\
\cline{2-2}\cline{4-5}
&  $\ne \pm r^2$ & & $\cC_2\times\cC_2$  & $\sqrt{-k}$\\
\hline
 \multirow{2}{*}{$16$} & $1,2$ & {$\cC_4$} & $\cC_2\times\cC_4$  & $\sqrt{2}$ \\
 \cline{2-5}
 &  $\ne 1,2 $ & {$\cC_2$} & $\cC_2\times\cC_2\,,\cC_4\,,\cC_4$  & $\sqrt{2},\sqrt{k},\sqrt{2k}$ \\
\hline
 $7$ &$-$  & {$\cC_2$} & $\cC_2\times\cC_2$  & $\sqrt{-7}$\\
 \hline
 $28$ & $-$ & {$\cC_2$} & $\cC_2\times\cC_2$  & $\sqrt{7} $\\
 \hline
 $8$ & $-$ & {$\cC_2$} & $\cC_2\times\cC_2$  & $\sqrt{2}$\\
 \hline
$11$ &$-$  & {$\cC_1$} & $-$  &$-$ \\
 \hline
$19$  &$-$  & {$\cC_1$} & $-$  & $-$\\
 \hline
 $43$  & $-$ & {$\cC_1$} & $-$  &$-$ \\
 \hline
 $67$ &  $-$ & {$\cC_1$} & $-$  & $-$\\
 \hline
 $163$  &$-$  & {$\cC_1$} & $-$  & $-$\\
 \hline
 \end{tabular}
\end{table}

{\em Notation:} Let $n$ denote a positive integer, we will denote by $\cC_n=\Z/n\Z$ the cyclic group of order $n$. Given an elliptic curve $E:y^2=x^3+Ax+B$, $A,B\in K$, and a number field $K$ we denote by $j(E)$ its $j$-invariant, by $\Delta(E)$ the discriminant of that short Weierstrass model, and by $E(K)_{\tors}$ the torsion subgroup of the Mordell-Weil group of $E$ over $K$.

\

\section{Preliminary results}
In this section we introduce some basic known facts that will be used in the proofs of the above theorems.
\subsection{Twists}\label{sec:twists} 
Let $E:y^2=x^3+Ax+B$ be an elliptic curve defined over $\Q$. Then any elliptic curve defined over $\Q$ isomorphic over $\Qbar$ to $E$ has a short Weierstrass model of the form:
$$
\begin{array}{lclcl}
\mbox{(i)} & & E^d:y^2=x^3+d^2Ax+d^3B & & \mbox{if $j(E)\ne 0,1728$,} \\
\mbox{(ii)} & & E^d:y^2=x^3+d Ax & & \mbox{if $j(E) =1728$,}\\
\mbox{(iii)} & & E^d:y^2=x^3+d B & & \mbox{if $j(E) = 0$,}
\end{array}
$$ 
where $d$ is an integer in $\Q^*/(\Q^*)^{n(E)}$ and $n(E)=2$ (resp. $4$ or $6$) if $j(E)\ne 0,1728$ (resp. $j(E)= 1728$ or $j(E)= 0$) (cf. \cite[X \S 5]{Silverman}). The elliptic curve $E^d$ is called the $d$-twist of $E$, and in the particular case $j(E)\ne 0,1728$ it is called the $d$-quadratic twist of $E$.

\subsection{Division polynomials} One of the main tools that we will use in this paper are the division polynomials of an elliptic curve (cf. \cite[\S 3.2]{Washington}). {Let $n$ be a positive integer and $E$ be an elliptic curve, we define the primitive $n$-division polynomial $\Psi_n(x)$ recursively, by dividing the (classical) $n$-division polynomial by the primitive $m$-division polynomial for all proper factors $m$ of $n$.} Then $\Psi_n(x)$ is characterized by the property that its roots are the $x$-coordinates of the points of exact order $n$ of $E$. Note that in general the $n$-division polynomial is defined so that its roots are the $x$-coordinates of the points of order dividing $n$, that is, the points in $E[n]$. 
\subsection{Quadratic twists}\label{sec:quadtwist}  Let $E:y^2=x^3+Ax+B$ be an elliptic curve defined over $\Q$, $d\in \Q$ squarefree, $E^{[d]}:y^2=x^3+Ad^2x+Bd^3$, and $E^{(d)}:dy^2=x^3+Ax+B$. Note that if $j(E)\ne 0,1728$ then $E^{[d]}=E^d$ is its $d$-quadratic twist, meanwhile if $j(E)=1728$ then $E^{[d]}=E^{d^2}$; and if $j(E)=0$ then $E^{[d]}=E^{d^3}$. We have the following isomorphisms:
$$
\begin{array}{rcl}
E & \longrightarrow & E^{(d)}\\
(x,y)& \mapsto &(x,y/\sqrt{d})
\end{array}
\qquad\mbox{and}\qquad
\begin{array}{rcl}
E^{(d)}& \longrightarrow & E^{[d]}\\
(x,y)& \mapsto &(dx,d^2y)\,.
\end{array}
$$

In the special case of quadratic twists there are two interesting results that will be useful in the sequel.
\begin{enumerate}
\item\label{qua1} The composition of the above two maps gives an isomorphism between $E$ and $E^{[d]}$ such that if $P=(\alpha,\beta)\in E$ then $P'=(d\alpha,d^{3/2}\beta)\in E^{[d]}$. In particular if $P\in E[n]$ then $P'\in E^{[d]}[n]$ for any positive integer $n$. On the other hand, suppose that we have $P=(\alpha,\beta)\in E[n]$ with $\alpha\in K$, in particular $\alpha$ is a root of $\Psi_n(x)$. In order to determine if there exist an square free integer $d$ such that $P'\in E^{[d]}(K)[n]$ we only need to check if there exists $\gamma\in K$ such that $\alpha^3+A\alpha+B=d\gamma^2$.
\item\label{qua2} If $n$ is an odd integer: $E(\Q(\sqrt{d}))[n]\simeq E(\Q)[n] \oplus E^{[d]}(\Q)[n]$.
\end{enumerate}

\subsection{Elliptic curves over $\Q$ with CM}\label{sec:CM}
Thanks to the classical theory of complex multiplication, there are only thirteen classes (up to $\Qbar$-isomorphism) of elliptic curves defined over $\Q$ with CM (cf. \cite[A \S 3]{SilvermanAd}). Each of these thirteen $j$-invariants corresponds to an order $R=\Z+\mathfrak{f}\,\mathcal{O}_K$ of conductor $\mathfrak{f}$ in a quadratic imaginary field $K=\Q(\sqrt{-D})$, where $ \mathcal{O}_K$ is the ring of integer of $K$. The thirteen possibilities are 
$$
(-D,\mathfrak{f})\in \left\{
\begin{array}{c}
(-3,1),(-3,2),(-3,3),(-4,1),(-4,2),(-7,1),(-7,2)\\
(-8,1),(-11,1),(-19,1),(-43,1),(-67,1),(-163,1)
\end{array}
\right\}.
$$
For the sake of simplicity we will denote by $\mathfrak{cm}$ the absolute value of the discriminant of the CM quadratic order $R$, that is $\mathfrak{cm}=D\cdot \mathfrak{f}^2$. Table \ref{isoCM} gives a representative elliptic curve $E_{\mathfrak{cm}}$ over $\Q$ for each $\mathfrak{cm}$.

\begin{table}[!h]
\caption{Isomorphism classes of elliptic curves defined over $\Q$ with CM.}\label{isoCM}
\begin{tabular}{|c|l|c|}
\hline
$\mathfrak{cm}$ &$E_{\mathfrak{cm}}\,:\,y^2=f_{\mathfrak{cm}}(x)$ &  $j(E_{\mathfrak{cm}})$ \\
\hline
$3$  & $\quad y^2=x^3 + 1$ & $0$\\
\hline
 $12$ & $\quad y^2=x^3-15 x+22$ & $2^4\cdot 3^3\cdot 5^3$ \\
\hline
$27$ & $\quad y^2=x^3 - 480 x + 4048$ & $-2^{15}\cdot  3\cdot  5^3$\\
\hline
$4$  & $\quad y^2=x^3+x$ & $2^6\cdot  3^3=1728$\\
\hline
$16$ & $\quad y^2=x^3 - 11 x + 14$ & $2^3\cdot  3^3\cdot  11^3$ \\
\hline
$7$  & $\quad y^2=x^3 - 2835x - 71442$ & $-3^3\cdot  5^3$\\
\hline
$28$ & $\quad y^2=x^3-595x+5586$ & $3^3\cdot  5^3\cdot  17^3$\\
\hline
$8$ &  $\quad y^2=x^3 - 4320 x + 96768$ & $2^6\cdot  5^3$\\
\hline
$11$ &  $\quad y^2=x^3 - 9504 x + 365904$ & $-2^{15}$\\
\hline
$19$ &  $\quad y^2=x^3 - 608 x + 5776 $ & $-2^{15}\cdot  3^3$\\
\hline
$43$ &  $\quad y^2=x^3 - 13760 x + 621264$ & $-2^{18}\cdot  3^3\cdot  5^3$\\
\hline
$67$ &  $\quad y^2=x^3 - 117920 x + 15585808 $ & $-2^{15}\cdot  3^3\cdot  5^3\cdot  11^3$\\
\hline
$163$  & $\quad y^2=x^3 - 34790720 x + 78984748304\quad $ & $-2^{18}\cdot  3^3\cdot  5^3\cdot  23^3\cdot  29^3$\\
\hline
\end{tabular}
\end{table}

Let $E$ be an elliptic curve defined over  $\Q$ with CM, by \S \ref{sec:twists} we have that $E$ is $\Q$-isomorphic to a curve $E^k_{\mathfrak{cm}}$ for some $\mathfrak{cm}$ as in Table \ref{isoCM}, and $k$ an integer in $\Q^*/(\Q^*)^{n(E)}$. Then $k$ and $\mathfrak{cm}$ are uniquely determined by $E$. We call them the CM-invariants of the elliptic curve $E$.

\section{Torsion over $\Q$}\label{sec:torsQ}
Let $E$ be an elliptic curve defined over $\Q$ with CM. In this section we compute the torsion subgroup of $E$ depending on its CM-invariants $(\mathfrak{cm},k)$. Note that this is a well-known\footnote{For example: the case $\mathfrak{cm}=3$ was first computed by Fueter  \cite{Fueter1930}; the case $\mathfrak{cm}=4$ in  \cite[\S 3]{Olson74}.} result but for the sake of completeness we include here the details of the proofs since they are going to be useful for the study of the torsion growth to quadratic fields. 

Suppose that $E$ has CM-invariants $(\mathfrak{cm},k)$, then $E$ is $\Q$-isomorphic to $E^k_\mathfrak{cm}$. Thanks to Olson's classification of $\Phi^{\cm}(1)$, in order to determine $E(\Q)_{\tors}$ we only need to study if the $2$-, $3$- and $4$-division polynomials have rational roots. Note that if the $n$-division polynomial of $E$ has no rational roots, then neither does the $n$-division polynomial of any quadratic twist of $E$. In the cases where $j(E)\notin\{0,1728\}$ there are only quadratic twists. In particular it is only necessary to study the $2$-, $3$- and $4$-division polynomials for $E_\mathfrak{cm}$. In the following cases the $n$-division polynomial $\Psi_n(x)$ refers to the elliptic curve $E_\mathfrak{cm}$.\\
\indent $\bullet$ $\mathfrak{cm}\in \{11,19,43,67,163\}$: $E(\Q)_{\tors}$ is trivial since $\Psi_2(x)$ and $\Psi_3(x)$ have no rational roots.\\
\indent $\bullet$ $\mathfrak{cm}\in \{7,28,8\}$: $E(\Q)_{\tors}\simeq \cC_2$ since $\Psi_2(x)$ has only one rational root and, $\Psi_3(x)$ and $\Psi_4(x)$ have no rational roots.\\
\indent $\bullet$  $\mathfrak{cm}=16$: $\Psi_3(x)$ has no rational roots, but $\Psi_2(x)$ has only one rational root. Let us check if there are points of order $4$. $\Psi_4(x)$ has two rational roots $r_1,r_2\in\Q$ and $f_{16}(r_i)=i s_i^2$ for $s_1,s_2\in\Q$. That is, only for $k=1,2$ the $k$-quadratic twist has points of order $4$. Therefore $E^k_{16}(\Q)_{\tors}\simeq \cC_2$ for $k\ne 1,2$ and $E^k_{16}(\Q)_{\tors}\simeq \cC_4$ for $k=1,2$.\\
\indent $\bullet$ $\mathfrak{cm}=27$: $\Psi_2(x)$ has no rational roots and $\Psi_3(x)$ has only one rational root $r\in \Q$. Now, $f_{27}(r)=s^2$ for some $s\in \Q$. Therefore, $E^k_{27}(\Q)_{\tors}\simeq \cC_3$ if $k=1$ and $E^k_{27}(\Q)_{\tors}$ is trivial if $k\ne 1$.\\
\indent $\bullet$ $\mathfrak{cm}=12$: $\Psi_2(x)$ has only one rational root and $\Psi_4(x)$ has not rational roots. Now, $\Psi_3(x)$ has only one rational root $r\in \Q$ and $f_{12}(r)=s^2$ for some $s\in\Q$. Therefore, $E^k_{12}(\Q)_{\tors}\simeq \cC_6$ if $k=1$ and $E^k_{12}(\Q)_{\tors}\simeq\cC_2$ if $k\ne 1$.\\

Let us study the non-quadratic twists:\\

\indent $\bullet$ $\mathfrak{cm}=4$. In this case the elliptic curve $E$ is $\Q$-isomorphic to $E_4^k:y^2=x^3+kx$ for some $k\in \Q^*/(\Q^*)^{4}$. The point $(0,0)\in E_4^k(\Q)$ is of order $2$ for any $k$. Let us see if there are points of order $3$: $\Psi_3(x)=3 x^4 + 6kx^2 - k^2$, then $z=-1/3(3\pm 2\sqrt{3})k$ are the roots of the polynomial $\Psi_3(\sqrt{x})$, but $z\ne x^2 $ for any $x,k\in \Q$. Therefore there are no points of order $3$ for any $k$. Now, let us check the existence of points of order $4$: $\Psi_4(x)=2(x^2 - k)(x^4 + 6kx^2 + k^2)$. Analogously to the previous case, the factor $x^4 + 6kx^2 + k^2$ has no rational roots for any $k$. But the first factor $x^2-k$ has rational roots if $k=r^2$ for some $r\in\Q$, in that case $x=\pm r$. Then $f_4(\pm r)=\pm 2 r^3$ is a rational square if and only if $r=\pm 2$. Therefore we conclude that $E^k_{4}(\Q)_{\tors}\simeq \cC_4$ if $k=4$; $E^k_{4}(\Q)_{\tors}\simeq \cC_2\times \cC_2$ if $k=-r^2$; and $E^k_{4}(\Q)_{\tors}\simeq \cC_2$ otherwise. \\
\indent $\bullet$ $\mathfrak{cm}=3$. In this case the elliptic curve $E$ is $\Q$-isomorphic to $E_3^k:y^2=x^3+k$ for some $k\in \Q^*/(\Q^*)^{6}$. It has points of order $2$ if and only if $k=s^3$ for some squarefree $s\in \Q$. In this case it is not possible to have full $2$-torsion over $\Q$ since $x^3+r^3=(x-r)(x^2+rx+r^2)$ and $x^2+rx+r^2$ is irreducible over $\Q$ for any $r\in\Q$. Let us study if there are points of order $3$. We look at the primitive $3$-division polynomial $\Psi_3(x)=3x(x^3 + 4k)$. If $x=0$ then $f_{3}(0)=k$. Therefore there is a rational point of order $3$ with $x$-coordinate $0$ if and only $k=r^2$ for some $r\in \Q$. If $x^3 + 4k =0$ then $k=\pm 2r^3$ for some squarefree $r\in\Q$ and $x= \mp 2 r$. Since $f_3(\mp 2 r)= \pm 6r^3$, a similar argument to the case $x=0$ allows us to conclude $r=\mp 6$. That is $k=-432$. We have obtained that there are points of order $3$ if and only if $k=r^2$ or $k=-432$. Therefore $E^k_{3}(\Q)_{\tors}\simeq \cC_6$ if $k=s^3$ and $k=r^2$, that is $k=1$; $E^k_{3}(\Q)_{\tors}\simeq \cC_3$ if $k=r^2\ne 1$ or $k=-432$;  $E^k_{3}(\Q)_{\tors}\simeq \cC_2$ if $k=r^3\ne 1$; and $E^k_{3}(\Q)_{\tors}$ is trivial otherwise. 

\

We have proved the following result:

\begin{proposition}\label{ProptorsQ}
Let $E$ be an elliptic curve defined over $\Q$ with CM. Table \ref{torsQ} gives an explicit description of the torsion subgroup of $E$ depending only on its CM-invariants $k$ and $\mathfrak{cm}$. 
\end{proposition}

\begin{table}[ht!]
\caption{Torsion of elliptic curves defined over $\Q$ with CM.}\label{torsQ}
\begin{tabular}{cccc}
\begin{tabular}{|c|c|c|}
\hline
$\mathfrak{cm}$ & $k$ & $E^k_{\mathfrak{cm}}(\Q)_{\tors}$  \\
\hline
 \multirow{4}{*}{$3$}  &  $1$ & $\cC_6$ \\
\cline{2-3}
& $-432,r^2\ne 1$&  $\cC_3$\\
\cline{2-3}
 &  $r^3 \ne 1$&  $\cC_2$\\
\cline{2-3}
 &  $\ne r^2,r^3,-432$ &  $\cC_1$ \\
\hline
 \multirow{2}{*}{$12$} & $1$ & $\cC_6$ \\
\cline{2-3}
 &   $\ne 1$ & $\cC_2$ \\
\hline
\multirow{2}{*}{$27$} & $1$ & $\cC_3$ \\
\cline{2-3}
 &   $\ne 1$ & $\cC_1$ \\
\hline
 \multirow{3}{*}{$4 $} & $4$ & $\cC_4$ \\
\cline{2-3}
&   $-r^2$ & $\cC_2\times \cC_2$ \\
\cline{2-3}
&   $\ne 4, -r^2$ & $\cC_2$ \\
\hline
\multirow{2}{*}{$16 $} & $1,2$ & $\cC_4$ \\
\cline{2-3}
&   $\ne 1,2$ & $\cC_2$ \\
\hline
\end{tabular}
& & & 
\begin{tabular}{|c|c|}
\hline
$\mathfrak{cm}$ & $E^k_{\mathfrak{cm}}(\Q)_{\tors}$  \\
\hline
$7 $ &  \multirow{3}{*}{$\cC_2$} \\
\cline{1-1}
  $28 $  &  \\
\cline{1-1}
 $ 8$ &  \\
\hline
$ 11$ & \multirow{5}{*}{$\cC_1$} \\
\cline{1-1}
 $19 $  &  \\
\cline{1-1}
 $43 $ &  \\
\cline{1-1}
 $67 $ & \\
\cline{1-1}
 $163 $ &  \\
\hline
\multicolumn{2}{c}{}\\
\multicolumn{2}{c}{}\\
\multicolumn{2}{c}{}\\
\multicolumn{2}{c}{}\\
\multicolumn{2}{c}{}\\
\end{tabular}
\end{tabular}
\end{table}

\section{Torsion growth over quadratic fields}

\subsection{Proof of Theorem \ref{teo1}} 
Let $H\in \Phi^{\cm}(2)\setminus\{\cC_5,\cC_7\}$. Table \ref{d2} shows examples of elliptic curves $E$ defined over $\Q$ with CM and quadratic fields $K$ such that $E(K)_{\tors}\simeq H$. Now, by Olson's classification we know that there are no elliptic curves with CM defined over $\Q$ with points of order $5$ (resp. $7$) over $\Q$. Finally, (\ref{qua2}) in \S\ref{sec:quadtwist} shows that there cannot be points of order $5$ (resp. $7$) over a quadratic field. Therefore $\cC_{10}, \cC_7\notin \Phi_{\Q}^{\cm}(2)$.

{
\begin{remark}
Let $K$ be a quadratic field, and let $E$ be an elliptic curve defined over $K$ with CM by a quadratic order of discriminant $-\mathfrak{cm}$ such that $E(K)_{\tors}\notin\left\{ \cC_1\,,\,  \cC_2\,,\,  \cC_3\,,\,  \cC_4\,,\,  \cC_6\,,\, \cC_2\times\cC_2\right\}$. Bourdon, Clark and Stankewicz \cite[Theorem 1.4]{BCS17} have proved that $K$ is the quadratic field listed below, and over that field $E$ is isomorphic to $\mathcal{E}_{\alpha,\beta}:y^2+(1-\alpha)xy-\beta y=x^3-\beta x^2$ where $\mathcal{E}_{0,0}:x^3+y^3=z^3$:
\begin{center}
    \begin{tabular}{ |c|c|c|c|c|c|}
    \hline
$K$ &   $\alpha$ & $\beta$ & $\mathfrak{cm}$ & $E(K)_{\tors}$ & Base Change from $\Q$? \\  \hline
    $\Q(\sqrt{-3})$ & $0$ & $0$ & $3$ & $\cC_{3} \times \cC_3$ & $E_3^{16}$ and $E_3^{-432}$ \\   [.5 ex]
     \hline
    $\Q(\sqrt{-1})$ & $-\frac{1}{8} $ & $0$ & $4$ & $\cC_2 \times \cC_4$ & $E_4^4$ and $E_4^{-1}$\\   [.5 ex]
    $\Q(\sqrt{2})$ & $1+\frac{3}{4}\sqrt{2}$ & $0$ & $4$ & $\cC_2 \times \cC_4$& $E_4^{-4}$ and $E_4^{-1}$ \\   [.5 ex]
    $\Q(\sqrt{2})$ & $-\frac{1}{32}$ & $0$ & $16$ & $\cC_2 \times \cC_4$ & $E_{16}^1$ and $E_{16}^2$ \\   [.5 ex]
    $\Q(\sqrt{2})$ & $\frac{1+\sqrt{2}}{8}$ & $0$ & $8$ & $\cC_2 \times \cC_4$ & No \\   [.5 ex]
    $\Q(\sqrt{-7})$ & $\frac{-31+3\sqrt{-7}}{512}$ & $0$ & $7$ & $\cC_2 \times \cC_4$ & No \\   [.5 ex]
    $\Q(\sqrt{-7})$ & $\frac{-1+3\sqrt{-7}}{32}$ & $0$ & $7$& $\cC_2 \times \cC_4$ & No \\   [.5 ex]
     \hline
    $\Q(\sqrt{-3})$ & $-\frac{2}{9}$ & $-\frac{1}{3}$ & $3$ & $\cC_2\times \cC_6$ & $E_3^{1}$ and $E_3^{-27}$ \\   [.5 ex]
    $\Q(\sqrt{3})$ & $\frac{1-\sqrt{3}}{9}$ & $\frac{-2+\sqrt{3}}{3}$ & $12$ & $\cC_2\times \cC_6$ & $E_{12}^{1}$ and $E_{12}^{3}$ \\   [.5 ex]
    $\Q(\sqrt{3})$ & $\frac{4}{9}$ & $\frac{1}{3}$ & $12$ & $\cC_2 \times \cC_6$ & $E_{12}^{1}$ and $E_{12}^{3}$ \\   [.5 ex]
     \hline
    $\Q(\sqrt{-3})$ &$ \frac{-1+\sqrt{-3}}{2}$ & $-1$ & $3$ & $\cC_7$ & No\\   [.5 ex]
     \hline
    $\Q(\sqrt{-1})$ & $\sqrt{-1}$ & $\sqrt{-1}$ & $4$ & $\cC_{10}$ & No \\   [.5 ex]
        \hline
 \end{tabular}
\end{center}

\

In the last column we show whether the elliptic curve $\mathcal{E}_{\alpha,\beta}$ is a base change of an elliptic curve over $\Q$ to the quadratic field $K$. In the affirmative case there always appear two elliptic curves defined over $\Q$. These elliptic curves are isomorphic over the quadratic field $K$.  Note that there is a typo in \cite[Theorem 1.4]{BCS17} since the two elliptic curves in the above table with $\mathfrak{cm}=12$ are isomorphic over $K=\Q(\sqrt{3})$, so there should appear only one. 
\end{remark}
}

\subsection{Proof of Theorems \ref{teo2} and \ref{teo3}} 
The first part of Theorem \ref{teo2} is to determine the set $\Phi^{\cm}_\Q(2,G)$. If $G\in \Phi^{\cm}(1)$, $G\ne  \cC_2\times\cC_2$, Table \ref{d2} shows examples of elliptic curves $E$ defined over $\Q$ with CM and quadratic fields $K$ for any possible torsion structure in $\Phi_\Q(2,G)\cap \Phi_\Q^{\cm}(2)$  (cf. \cite[Theorem 2]{GJT14}). If $G=\cC_2\times\cC_2$, $\Phi_\Q(2,G)\cap \Phi_\Q^{\cm}(2)=\{\cC_2\times\cC_2,\cC_2\times\cC_4,\cC_2\times\cC_6 \}$  (cf. \cite[Theorem 2]{GJT14}). Then to finish the first part of Theorem \ref{teo2} we must prove that if $E$ is an elliptic curve defined over $\Q$ with CM such that $E(\Q)_{\tors}\simeq \cC_2\times\cC_2$ then the torsion cannot grow to $\cC_2\times\cC_6$ over a quadratic field. By Proposition \ref{ProptorsQ} (see Table \ref{torsQ}), $E$ should have  $\fcm=4$. But since there are no elliptic curves in the family $E^k_4$ with points of order $3$ over $\Q$ there cannot be points of order $3$ over a quadratic field (by (\ref{qua2}) in \S\ref{sec:quadtwist}).

The second part of Theorem \ref{teo2} is to determine $\mathcal{H}^{\cm}_\Q(2,G)$ for any $G\in \Phi^{\cm}(1)$. Notice that this is a direct consequence of Theorem \ref{teo3}, then we will prove Theorem \ref{teo3}  first. That is, we are going to justify all the entries in Table \ref{d2} following a similar argument to the one in section \S \ref{sec:torsQ}. 

For any elliptic curve $E$ defined over $\Q$ with CM, Proposition \ref{ProptorsQ} gives an explicit description of $G=E(\Q)_{\tors}$ in terms of its CM-invariants. Now thanks to the classification of $\Phi^{\cm}_\Q(2,G)$  we know the possible torsion growth over quadratic fields. In this case we only need to study if the $2$-, $3$- and $4$-division polynomials have linear or quadratic factors.

Remember that if $E$ has CM-invariants $(\mathfrak{cm},k)$, then $E$ is $\Q$-isomorphic to $E^k_\mathfrak{cm}$ and in the cases where $\mathfrak{cm}\notin\{3,4\}$ there are only quadratic twists. In particular it is only necessary to study the $2$-, $3$- and $4$-division polynomials for $E_\mathfrak{cm}$. In the following cases $\Psi_n(x)$ denotes the $n$-division polynomial of $E_\mathfrak{cm}$.

\indent $\bullet$ $\mathfrak{cm}\in \{11,19,43,67,163\}$:  $E^k_{\fcm}(\Q)_{\tors}\simeq \cC_1$ and, since $\Phi^{\cm}_\Q(2,\cC_1)=\{\cC_1,\cC_3\}$, there is no torsion growth over quadratic fields.\\
\indent $\bullet$ $\mathfrak{cm}\in \{7,28,8\}$: $E^k_{\fcm}(\Q)_{\tors}\simeq \cC_2$ and $\Phi^{\cm}_\Q(2,\cC_2)=\{\cC_2,\cC_4,\cC_6,\cC_2\times\cC_2,\cC_2\times\cC_6\}$. Similarly to the previous case there cannot be points of order $3$ over quadratic fields, that is, neither torsion growth to $\cC_6$  nor $\cC_2\times\cC_6$. Now since the torsion over $\Q$ is $\cC_2$ we have that the full $2$-torsion is defined over $\Q(\sqrt{\Delta(E_{\fcm})})$. In our cases we have $\Q(\sqrt{\Delta({E_{7}})})=\Q(\sqrt{-7})$, $\Q(\sqrt{\Delta({E_{28})}})=\Q(\sqrt{7})$, and $\Q(\sqrt{\Delta({E_{8}})})=\Q(\sqrt{2})$. Finally let us check if there is torsion growth to $\cC_4$.
\begin{itemize}
\item $\fcm=7$:  $\Psi_4(x)=2(x^2 - 126x - 5103)(x^2 + 567)(x^2 + 126x + 6237)$. The second and third quadratic irreducible factors have squarefree part of the discriminant equal to $-7$. Then a possible point of order $4$ should be defined over the field of definition of the full $2$-torsion. But this is impossible since $\cC_2\times \cC_4$ is not a subgroup of a group in $\Phi^{\cm}_\Q(2,\cC_2)$. Now, $\alpha=63+36\sqrt{7}$ is a root of the first quadratic factor of $\Psi_4(x)$. We have $f_7(\alpha)=\sqrt{7}u( 2^2 3^3\sqrt{7})^2$, where $u=8+3\sqrt{7}$ is a fundamental unit of the quadratic field $\Q(\sqrt{7})$.  Therefore $f_7(\alpha)\ne d \beta^2$ for any $d\in \Q$ and $\beta\in \Q(\sqrt{7})$. This proves that there are no points of order $4$ over any quadratic field.
\item $\fcm=28$:  $\Psi_4(x)$ has only one irreducible factor of degree $\le 2$. One of its roots is $\alpha=14+\sqrt{-7}$ and we have $f_{28}(\alpha)=-7\sqrt{-7}(1+\sqrt{-7})^4/4$. Similarly to the previous case: $f_{28}(\alpha)\ne d \beta^2$ for any $d\in \Q$ and $\beta\in \Q(\sqrt{-7})$ and there are no points of order $4$ over quadratic fields.
\item $\fcm=8$:  There is only one irreducible factor of $\Psi_4(x)$ of degree $\le 2$. In this case its roots are defined over $\Q(\sqrt{2})$. Since $\Q(\sqrt{\Delta({E_{8}})})=\Q(\sqrt{2})$ we obtain that there are no points of order $4$.
\end{itemize}
We have proved that there is only torsion growth to $\cC_2\times \cC_2$ over $\Q(\sqrt{\Delta({E_{\fcm}})})$.\\
\indent $\bullet$  $\mathfrak{cm}=16$: If $k=1,2$ then $E^k_{16}(\Q)_{\tors}\simeq \cC_4$. Fot these cases the torsion only grows to $\cC_2\times\cC_4$ over $\Q(\sqrt{\Delta({E_{16}})})=\Q(\sqrt{2})$ since $\Phi^{\cm}_\Q(2,\cC_4)=\{\cC_4,\cC_2\times\cC_4\}$. Now if $k\ne 1,2$, then $E^k_{16}(\Q)_{\tors}\simeq \cC_2$. By (\ref{qua2}) at section \S\ref{sec:quadtwist}) there are not points of order $3$ over quadratic fields since there are not points of order $3$ over $\Q$ for any quadratic twist of $E_{16}$. Finally, let us study if there are points of order $4$ over some quadratic field. The factorization of the $4$-division polynomial in irreducible factors is: $\Psi_4(x)=2(x-1)(x-3)(x^4 + 4x^3 - 42x^2 + 100x - 79)$. Now $f_{16}(1)=4$ and $f_{16}(3)=8$. Therefore there are points of order $4$ over the quadratic fields $\Q(\sqrt{k})$ and $\Q(\sqrt{2k})$. Since $k\ne 1,2$, the torsion subgroup over those quadratic fields is isomorphic to $\cC_4$.\\
\indent $\bullet$ $\mathfrak{cm}=27$: If $k=1$, then $E^1_{27}(\Q)_{\tors}\simeq \cC_3$. Since $\Phi^{\cm}_\Q(2,\cC_3)=\{\cC_3,\cC_3\times\cC_3\}$ the torsion can only grow to $\cC_3\times\cC_3$. But this can only happen over $\Q(\sqrt{-3})$. We compute that the torsion subgroup over $\Q(\sqrt{-3})$ is $\cC_3$ too. Then there is no torsion growth for $k=1$. Now suppose $k\ne 1$, then $E^k_{27}(\Q)_{\tors}\simeq \cC_1$ and $\Phi^{\cm}_\Q(2,\cC_1)=\{\cC_1,\cC_3\}$. By (\ref{qua2}) in \S\ref{sec:quadtwist} with $n=3$ we have that the torsion growth to $\cC_3$ over a quadratic field only over $\Q(\sqrt{k})$.\\
\indent $\bullet$ $\mathfrak{cm}=12$: If $k=1$, then $E^1_{12}(\Q)_{\tors}\simeq \cC_6$ and $\Phi^{\cm}_\Q(2,\cC_6)=\{\cC_6,\cC_2\times\cC_6\}$. Thus the torsion only grows to $\cC_2\times\cC_6$ over $\Q(\sqrt{\Delta({E^1_{12}}}))=\Q(\sqrt{3})$. Now if $k\ne 1$, then $E^k_{12}(\Q)_{\tors}\simeq \cC_2$ and $\Phi^{\cm}_\Q(2,\cC_2)=\{\cC_2,\cC_4,\cC_6,\cC_2\times\cC_2,\cC_2\times\cC_6\}$. We have that $\cC_2\times\cC_2$ is isomorphic to a subgroup of the torsion subgroup over $\Q(\sqrt{\Delta({E^1_{12}})})=\Q(\sqrt{3})$. A similar argument to the case $\mathfrak{cm}=27$ allows us to determine that there are points of order $3$ over $\Q(\sqrt{k})$. Therefore if $k=3$, the torsion grows only to $ \cC_2\times\cC_6$ over $\Q(\sqrt{3})$. Finally, if $k\ne 1,3$ we have torsion growth to $\cC_2\times\cC_2$  over $\Q(\sqrt{3})$, and $\cC_6$  over $\Q(\sqrt{k})$. It remains to check that there are no points of order $4$ over quadratic fields. $\Psi_4(x)$ has only one irreducible factor of degree $\le 2$. One of its roots is $\alpha=2+\sqrt{-3}$ and we have $f_{12}(\alpha)=-\sqrt{-3}(3-\sqrt{-3})^2$. Therefore: $f_{12}(\alpha)\ne d \beta^2$ for any $d\in \Q$ and $\beta\in \Q(\sqrt{-3})$. This proves that there are no points of order $4$ over any quadratic fields.

\

Finally we deal with the non-quadratic twists:

\

\indent $\bullet$ $\mathfrak{cm}=4$. We split the proof depending on the torsion over $\Q$:
\begin{itemize}
\item $E_4^{k}(\Q)_{\tors}\simeq \cC_2\times\cC_2$ if $k=-r^2$. We have $\Phi^{\cm}_\Q(2,\cC_2\times\cC_2)=\{\cC_2\times\cC_2,\cC_2\times\cC_4\}$. Let us study if there are points of order $4$ over a quadratic field. Note that, since $E_4^{-r^2}$ is the $r$-quadratic twist of $E_4^{-1}$, it is enough to study the factorization of the $4$-division polynomial $\Psi_4(x)$ of $E_4^{-1}$. The polynomial $\Psi_4(x)$ has the roots $\pm i$, $\pm 1\pm \sqrt{2}$ and evaluating the polynomial $x^3-x$ for these values we obtain:
$$
\qquad i^3-i=(i-1)^2\,,\quad (1+\sqrt{2})^3-(1+\sqrt{2})=(2+\sqrt{2})^2 \,,\quad (-1+\sqrt{2})^3-(-1+\sqrt{2})=-(2-\sqrt{2})^2.
$$
Therefore there are points of order $4$ over a quadratic fields if and only if $r=\pm 1$ over $\Q(i),\Q(\sqrt{2})$ or $r=\pm 2$ over $\Q(\sqrt{2})$. That is, over those quadratic fields and the corresponding values of $r$ we obtain that the torsion subgroup is isomorphic to $ \cC_2\times\cC_4$. For the rest of the values of $r$ there is no torsion growth over any quadratic field for the elliptic curve  $E_4^{-r^2}$. 
\item $E_4^{k}(\Q)_{\tors}\simeq \cC_4$ if $k=4$. There is only one possibility to grow over a quadratic field: $\cC_2\times\cC_4$ over $\Q(i)$, since $\Delta({E^4_{4}})=-2^{12}$. 
\item $E_4^{k}(\Q)_{\tors}\simeq \cC_2$ if $k\ne 4,-r^2$. By (\ref{qua2}) in \S\ref{sec:quadtwist} there are no points of order $3$ over quadratic fields since there are no points of order $3$ over $\Q$ for any value of $k$. First suppose $k=r^2$. Then $E_4^{r^2}$ is the $r$-quadratic twist of $E_4$. Let us study if there are points of order $4$ using the $4$-division polynomial of $E_4$: $\Psi_4(x)=2(x^2-1)(x^4+6x^2+1)$. Since $f_4(\pm 1)=\pm 2$ there are points of order $4$ over a quadratic field only in the case $\Q(\sqrt{\pm 2r})$. The last possibility for torsion growth is $\cC_2\times \cC_2$ over $\Q(\sqrt{-1})$. This finishes the case $k=r^2$.  Finally we deal with the general case:  $k\ne  r^2$. We have  $\cC_2\times \cC_2$ is isomorphic to a subgroup of the torsion subgroup over $\Q(\sqrt{-k})$ since $\Delta(E_4^k)=-k(8k)^2$. To finish the proof of this case we are going to prove that there are no points of order $4$ over quadratic fields. We have $\Psi_4(x)=2(x^2 - k)(x^4 + 6kx^2 + k^2)$. Let us denote by $g(x)$ the second factor, then $z=(-3\pm2 \sqrt{2})k$ are the roots of the polynomial $g(\sqrt{x})$, but $z\ne x^2 $ for any $x\in\Q(\sqrt{2})$ and $k\in \Q$. The first factor have the roots $x=\pm \sqrt{k}$, but $(\pm \sqrt{k})^3+k(\pm \sqrt{k})=\pm 2\sqrt{k}^3$ is never an square over $\Q(\sqrt{k})$. We conclude that there is only torsion growth over quadratic fields to $\cC_2\times \cC_2$ over $\Q(\sqrt{-k})$.\\
\end{itemize}

\indent $\bullet$ $\mathfrak{cm}=3$. Note that this case has been dealt by Dey\footnote{Note that there is a typo in Theorem 1(3) \cite{Dey19} since it is necessary to add the case $c=-27$ and $d\ne -3$ (in Dey's notation).} \cite{Dey19} with a slightly different approach. We split the proof depending on the torsion over $\Q$:
\begin{itemize}
\item $E_3^{k}(\Q)_{\tors}\simeq \cC_6$ if $k=1$. The torsion only grows to $\cC_2\times\cC_6$ over the quadratic field $\Q(\sqrt{\Delta(E_3^1)})=\Q(\sqrt{-3})$.
\item $E_3^{k}(\Q)_{\tors}\simeq \cC_3$ if $k=-432$ or $k=r^2\ne 1$. Here the torsion can only grow to $\cC_3\times\cC_3$  over $\Q(\sqrt{-3})$. We check that $E_3^{-432}(\Q(\sqrt{-3}))_{\tors}\simeq \cC_3\times \cC_3$. Now suppose $k=r^2\ne 1$. We must have all the roots of $\Psi_3(x)=3x(x^3+4r^2)$ defined over $\Q(\sqrt{-3})$. Therefore $r=4s^3$, but since $k\in\Q^*/(\Q^*)^6$ the unique possibility is $r=4$, i.e. $k=16$. We check that $E_3^{16}(\Q(\sqrt{-3}))_{\tors}\simeq \cC_3\times \cC_3$. For the rest of the values the torsion does not grow over quadratic fields.
\item $E_3^{k}(\Q)_{\tors}\simeq \cC_2$ if $k=r^3\ne 1$. Note that in this case $E_3^{k}$ is the $r$-quadratic twist of $E_3$, therefore it is enough to study the $n$-division polynomials $\Psi_n(x)$ of $E_3$. Since the torsion over $\Q$ is isomorphic to $\cC_2$, we could only have torsion growth $\cC_4,\cC_6,\cC_2\times \cC_2,\cC_2\times\cC_6$. Let us check if there are points of order $4$ over quadratic fields: The unique factor of $\Psi_4(x)$ of degree $\le 2$ is $g(x)=x^2 + 2x - 2$. We have that $\alpha= \sqrt{3}-1$ is a root of $g(x)$. Then $f_3(\alpha)=3\sqrt{3}u$, where $u=2-\sqrt{3}$ is a fundamental unit of the quadratic field $\Q(\sqrt{3})$. Therefore $f_3(\alpha)\ne r \beta^2$ for any $r\in \Q$ and $\beta\in \Q(\sqrt{3})$. In conclusion, there are no points of order $4$ over quadratic fields. Now we study if there are points of order $3$. We have $\Psi_3(x)=4x(x^3+4)$ and $f_3(0)=1$. Therefore there are points of order $3$ over $\Q(\sqrt{r})$. Finally we have full $2$-torsion over $\Q(\Delta(E_3))=\Q(\sqrt{-3})$. We conclude that there are torsion growth to $\cC_2\times\cC_2$ and $\cC_6$ if $k\ne -3$; and $\cC_2\times\cC_6$ if $r=-3$.
\item $E_3^{k}(\Q)_{\tors}\simeq \cC_1$ if  $k\ne r^2,r^3,-432$. We need to check only if there are points of order $3$ over quadratic fields. We have $\Psi_3(x)=3x(x^3+4k)$. If $x=0$ then $y^2=k$. Thus, the torsion grows to $\cC_3$ over $\Q(\sqrt{k})$. The second factor $x^3+4k$ has roots over a quadratic field if and only if $k=2 r^3$. In that case the root is $\alpha=-2r$ and we have $\alpha^3+k=-6r^3$ is a square over a quadratic field only over $\Q(\sqrt{-6r})$. Since $k=2r^3$ we must have $r\ne 2,-6$.
\end{itemize}

\

\begin{remark}
All the computations were done using \texttt{Magma} \cite{magma} and the source code is available at the author's webpage \cite{MagmaCode}.
\end{remark}

\

{\bf Acknowledgements.} The author would like to thank Harris B. Daniels, who read the earlier versions of this paper carefully. Finally, the author thanks the anonymous referees for their useful comments and suggestions.

\end{document}